# SECOND-ORDER ADJOINT SENSITIVITY ANALYSIS METHODOLOGY (2$^{nd}$-ASAM) FOR LARGE-SCALE NONLINEAR SYSTEMS: I. THEORY

Dan G. Cacuci

Department of Mechanical Engineering, University of South Carolina

E-mail: cacuci@cec.sc.edu

Corresponding author:

Department of Mechanical Engineering, University of South Carolina

300 Main Street, Columbia, SC 29208, USA

Email: cacuci@cec.sc.edu; Phone: (919) 909 9624;

**ABSTRACT**

This work presents the *second-order adjoint sensitivity analysis methodology* (*2$^{nd}$-ASAM*) *for nonlinear systems*, which yields exactly and efficiently the second-order functional derivatives of physical (engineering, biological, etc.) system responses (i.e., "system performance parameters") to the system's model parameters. The definition of "system parameters" used in this work includes all computational input data, correlations, initial and/or boundary conditions, etc. For a physical system comprising $N_\alpha$ parameters responses, forwards methods require a total of $\left( N_\alpha^2 / 2 + 3 N_\alpha / 2 \right)$ large-scale computations for obtaining all of the first- and second-order sensitivities, for all system responses. In contradistinction, the *2$^{nd}$-ASAM* requires one large-scale computation using the *first-level adjoint sensitivity system* (*1$^{st}$-LASS*) for obtaining all of the first-order sensitivities, followed by at most $N_\alpha$ large-scale computations using the *second-level adjoint sensitivity systems* (*2$^{nd}$-LASS*), for obtaining exactly all of the second-order sensitivities of a functional-type response. The construction, implementation and solution of the *2$^{nd}$-ASAM* requires very little additional effort beyond the construction of the *first-level adjoint sensitivity system* (*1$^{st}$-LASS*) needed for computing the first-order sensitivities. Furthermore, due to the symmetry properties of the second-order sensitivities, the *2$^{nd}$- ASAM* comprises the inherent automatic "solution verification" of the correctness and accuracy of the 2$^{nd}$-level adjoint functions used for the efficient and exact computation of the second-order sensitivities. The use of the *2$^{nd}$-ASAM* to compute exactly all of the second-order response sensitivities to model input parameters is expected to enable significant advances in related scientific disciplines, particularly the areas of uncertainty quantification and predictive modeling, including model validation, reduced-order modeling, data assimilation, model calibration and extrapolation.

**KEYWORDS:** second-order adjoint sensitivity analysis methodology (*2$^{nd}$-ASAM*); exact and efficient computation of first- and second-order functional derivatives; large-scale nonlinear systems.

## 1. INTRODUCTION

In a recent work [1], Cacuci (2015) has introduced the "**S**econd-**O**rder **A**djoint **S**ensitivity **A**nalysis **M**ethodology ($2^{nd}$-ASAM) *for **linear** systems*" for computing exactly and most efficiently the second-order functional derivatives of physical (engineering, biological, etc.) system responses (i.e., "system performance parameters") to the system's parameters. The term "system parameter" refers, in the most comprehensive sense, to all input data, correlations, initial and/or boundary conditions, etc. The "*$2^{nd}$-ASAM for linear systems*" considers nonlinear responses associated with physical systems modeled mathematically by systems of linear operator equations. The comparative discussion presented in [1] regarding the basic properties of the leading methods (deterministic and/or statistical) used for computing second-order sensitivities [2-8] and the fundamentally new and distinctive features of *$2^{nd}$-ASAM for linear systems* introduced in [1] highlighted the unparalleled efficiency of the *$2^{nd}$-ASAM for linear systems* for computing $2^{nd}$-order sensitivities *exactly*. Since the comparative discussion presented in [1] continues to remain valid in the context of the new *$2^{nd}$-ASAM for nonlinear systems* which will be introduced in this work, that discussion will not be repeated here. The efficiency of the "*$2^{nd}$-ASAM for linear systems*" for computing exactly first- and second-order sensitivities (i.e., functional derivatives) of model responses to model parameters has since been also demonstrated [9-11] in several recent applications to particle diffusion and heat transport problems.

Extending the work in [1], the present work introduces the "**S**econd-**O**rder **A**djoint **S**ensitivity **A**nalysis **M**ethodology ($2^{nd}$-ASAM) *for **nonlinear** systems*," which is a new method for computing exactly and efficiently *second-order functional derivatives of noninear system responses* (i.e., "system performance parameters" in physical, engineering, biological systems) *to the system's parameters characterizing large-scale nonlinear systems*. Just as in [1], the "*$2^{nd}$-ASAM for nonlinear systems*" builds on the *first-order adjoint sensitivity analysis methodology ($1^{st}$-ASAM) for linear and nonlinear systems* originally introduced by Cacuci in [12, 13], and further discussed in [14-16].

As is well known, response sensitivities to model parameters are needed in many applications, including:

(i) understanding the system by identifying and ranking the importance of model parameters in influencing the response under consideration;

(ii) determining the effects of parameter variations on the system's behavior;

(iii) improving the system design, possibly reducing conservatism and redundancy;

(iv) prioritizing possible improvements for the system under consideration;

(v) quantifying uncertainties in responses due to quantified parameter uncertainties (e.g., by using the method of "propagation of uncertainties" (see, e.g., [14]);

(vi) performing "predictive modeling" [17-24], including data assimilation, model calibration and extrapolation, for the purpose of obtaining best-estimate predicted results with reduced predicted uncertainties.

The work presented in this article is organized as follows: Section 2 briefly recalls the first-order adjoint sensitivity analysis methodology ($1^{st}$-ASAM), which was introduced in [12, 13] and which provides the foundation for the new $2^{nd}$-ASAM *for nonlinear systems*, which is presented in Section 3. As shown in Section 3, the $2^{nd}$-ASAM *for nonlinear systems* computes exactly and most efficiently all of the $2^{nd}$-order response sensitivities in at most $N_\alpha$ large-scale computations, as opposed to computing inexactly (e.g., via finite-differences) these $2^{nd}$-order response sensitivities in $\left( N_\alpha^2 / 2 + 3 N_\alpha / 2 \right)$ large-scale computations, as would be required by forward methods. Section 4 concludes this work by highlighting the significance of the $2^{nd}$-ASAM and noting that the sequel to this work [25] presents an illustrative paradigm application of the $2^{nd}$-ASAM to a nonlinear heat conduction problem.

## 2. BACKGROUND: THE FIRST-ORDER ADJOINT SENSITIVITY ANALYSIS METHODOLOGY ($1^{st}$-ASAM) FOR LARGE-SCALE NONLINEAR SYSTEMS

Consider that the physical system is represented mathematically by means of $N_u$ coupled *nonlinear* operator equations of the form

$$\mathbf{N}[\mathbf{u}(\mathbf{x}), \boldsymbol{\alpha}(\mathbf{x})] = \mathbf{Q}[\boldsymbol{\alpha}(\mathbf{x})], \quad \mathbf{x} \in \Omega_x, \tag{1}$$

where:

1. $\mathbf{x} = (x_1, \ldots, x_{J_x})$ denotes the $J_x$-dimensional phase-space position vector for the primary system; $\mathbf{x} \in \Omega_x \subset \mathbb{R}^{J_x}$, where $\Omega_x$ is a subset of the $J_x$-dimensional real vector space $\mathbb{R}^{J_x}$;

2. $\mathbf{u}(\mathbf{x}) = [u_1(\mathbf{x}), \ldots, u_{N_u}(\mathbf{x})]$ denotes a $N_u$-dimensional column vector whose components are the system's dependent (i.e., state) variables; $\mathbf{u}(\mathbf{x}) \in \mathcal{E}_u$, where $\mathcal{E}_u$ is a normed linear space over the scalar field $\mathcal{F}$ of real numbers;

3. $\boldsymbol{\alpha}(\mathbf{x}) = [\alpha_1(\mathbf{x}), \ldots, \alpha_{N_\alpha}(\mathbf{x})]$ denotes an $N_\alpha$-dimensional column vector whose components are the system's parameters; $\boldsymbol{\alpha} \in \mathcal{E}_\alpha$, where $\mathcal{E}_\alpha$ is also a normed linear space;

4. $\mathbf{Q}[\boldsymbol{\alpha}(\mathbf{x})] = [Q_1(\boldsymbol{\alpha}), \ldots, Q_{N_u}(\boldsymbol{\alpha})]$ denotes a $N_u$-dimensional column vector whose elements represent inhomogeneous source terms that depend either linearly or nonlinearly on $\boldsymbol{\alpha}$; $\mathbf{Q} \in \mathcal{E}_Q$, where $\mathcal{E}_Q$ is also a normed linear space; the components of $\mathbf{Q}$ may be operators, rather than just functions, acting on $\mathbf{u}(\mathbf{x})$, $\boldsymbol{\alpha}(\mathbf{x})$ and $\mathbf{x}$;

5. $\mathbf{N}[\mathbf{u}(\mathbf{x}), \boldsymbol{\alpha}(\mathbf{x})] \equiv [N_1(\mathbf{u}, \boldsymbol{\alpha}), \ldots, N_{N_u}(\mathbf{u}, \boldsymbol{\alpha})]$ denotes a $N_u$-component column vector whose components are *operators* (including differential, difference, integral, distributions, and/or infinite matrices) acting *nonlinearly* on $\mathbf{u}$ and $\boldsymbol{\alpha}$. For notational convenience, all vectors in this work are considered to be column vectors; transposition will be indicated by a dagger $(\dagger)$.

6. All of the equalities in this work are considered to hold in the weak ("distributional") sense, since the right-sides ("sources") of the various eqautions, including Eq. (1) may contain distributions ("generalized functions/functionals"), particularly Dirac-distributions and derivatives and/or integrals thereof.

In view of the definitions given above, $\mathbf{N}(\mathbf{u},\boldsymbol{\alpha})$ represents the mapping $\mathbf{N}:\mathcal{D}\subset\mathcal{E}\to\mathcal{E}_Q$, where $\mathcal{D}=\mathcal{D}_u\times\mathcal{D}_\alpha$, $\mathcal{D}_u\subset\mathcal{E}_u$, $\mathcal{D}_\alpha\subset\mathcal{E}_\alpha$, and $\mathcal{E}=\mathcal{E}_u\times\mathcal{E}_\alpha$. Note that an arbitrary element $\mathbf{e}\in\mathcal{E}$ is of the form $\mathbf{e}=(\mathbf{u},\boldsymbol{\alpha})$. If differential operators appear in Eq. (1), then a corresponding set of boundary and/or initial conditions (which are essential to define $\mathcal{D}$) must also be given; these boundary and/or initial conditions are represented in operator form as

$$\left[\mathbf{B}(\mathbf{u},\boldsymbol{\alpha})-\mathbf{A}(\boldsymbol{\alpha})\right]_{\partial\Omega_x}=\mathbf{0}, \quad \mathbf{x}\in\partial\Omega_x, \tag{2}$$

where $\partial\Omega_x$ denotes the boundary of $\Omega_x$, the operator $\mathbf{B}(\mathbf{u},\boldsymbol{\alpha})$ acts nonlinearly on both $\mathbf{u}$ and on the model parameters $\boldsymbol{\alpha}$, while $\mathbf{A}(\boldsymbol{\alpha})$ denotes an operator that acts nonlinearly on $\boldsymbol{\alpha}$.

The vector-valued function $\mathbf{u}(\mathbf{x})$ is considered to be the unique nontrivial solution of the physical problem described by Eqs. (1) and (2). The system response (i.e., result of interest), associated with the problem modeled by Eqs. (1) and (2) will be denoted here as $R(\mathbf{u},\boldsymbol{\alpha})$; in this work, $R(\mathbf{u},\boldsymbol{\alpha})$ is considered to be a real-valued *nonlinear functional* of $(\mathbf{u},\boldsymbol{\alpha})$, which can be generally represented in operator form as

$$R(\mathbf{u},\boldsymbol{\alpha}):\mathcal{D}_R\subset\mathcal{E}\to\mathcal{F}, \tag{3}$$

where $\mathcal{F}$ denotes the field of real scalars.

The nominal parameter values $\boldsymbol{\alpha}^0(\mathbf{x})$ are used in Eqs. (1) and (2) to obtain the nominal solution $\mathbf{u}^0(\mathbf{x})$ by solving these equations; mathematically, therefore, the nominal value $\mathbf{u}^0(\mathbf{x})$ of the state-function is obtained by solving

$$\mathbf{N}\left[\mathbf{u}^0(\mathbf{x}),\boldsymbol{\alpha}^0(\mathbf{x})\right]=\mathbf{Q}\left[\boldsymbol{\alpha}^0(\mathbf{x})\right], \quad \mathbf{x}\in\Omega_x, \tag{4}$$

$$\mathbf{B}(\mathbf{u}^0,\boldsymbol{\alpha}^0)=\mathbf{A}(\boldsymbol{\alpha}^0), \quad \mathbf{x}\in\partial\Omega_x. \tag{5}$$

Equations (4) and (5) represent the "*base-case*" or *nominal* state of the physical system; the superscript "zero" will be used in this work to denote "nominal" values. After solving Eqs. (4) and (5), the *nominal solution*, $\mathbf{u}^0(\mathbf{x})$, thus obtained is used to obtain the *nominal value*, $R(\mathbf{e}^0)$, of the response $R(\mathbf{e})$, using the nominal values $\mathbf{e}^0 = (\mathbf{u}^0, \boldsymbol{\alpha}^0) \in \mathscr{E}$ of the model's state function and parameters.

Consider now (a vector of) arbitrary variations $\mathbf{h} \equiv (\mathbf{h}_u, \mathbf{h}_\alpha) \in \mathscr{E} = \mathscr{E}_u \times \mathscr{E}_\alpha$, with $\mathbf{h}_u \equiv (\delta u_1, \ldots, \delta u_{N_u}) \in \mathscr{E}_u$ and $\mathbf{h}_\alpha \equiv (\delta \alpha_1, \ldots, \delta \alpha_{N_\alpha}) \in \mathscr{E}_\alpha$, around $\mathbf{e}^0 = (\mathbf{u}^0, \boldsymbol{\alpha}^0) \in \mathscr{E}$. The variation (sensitivity) of the response $R$ to variations $\mathbf{h}$ in the system parameters is given by the Gâteaux- (G)-differential $\delta R(\mathbf{e}^0; \mathbf{h})$ of the response $R(\mathbf{e})$ at $\mathbf{e}^0 = (\mathbf{u}^0, \boldsymbol{\alpha}^0)$ with increment $\mathbf{h}$, which is defined as

$$\delta R(\mathbf{e}^0; \mathbf{h}) \equiv \left\{ \frac{d}{d\varepsilon} \left[ R(\mathbf{e}^0 + \varepsilon \mathbf{h}) \right] \right\}_{\varepsilon=0} = \lim_{\varepsilon \to 0} \frac{R(\mathbf{e}^0 + \varepsilon \mathbf{h}) - R(\mathbf{e}^0)}{\varepsilon}, \qquad (6)$$

for $\varepsilon \in \mathscr{F}$, and all (i.e., arbitrary) vectors $\mathbf{h} \in \mathscr{E}$. When the response $R(\mathbf{e})$ is functional of the form $R: \mathscr{D}_R \to \mathscr{F}$, the sensitivity $\delta R(\mathbf{e}^0; \mathbf{h})$ is also an operator, defined on the same domain, and with the same range as $R(\mathbf{e})$. The G-differential $\delta R(\mathbf{e}^0; \mathbf{h})$ is related to the total variation $\left[ R(\mathbf{e}^0 + \varepsilon \mathbf{h}) - R(\mathbf{e}^0) \right]$ of $R$ at $\mathbf{e}^0$ through the relation $R(\mathbf{e}^0 + \varepsilon \mathbf{h}) - R(\mathbf{e}^0) = \delta R(\mathbf{e}^0; \mathbf{h}) + \Delta(\mathbf{h})$, with $\lim_{\varepsilon \to 0} \left[ \Delta(\varepsilon \mathbf{h}) \right]/\varepsilon = 0$.

As discussed in [1] and [2], *the most general definition of the first-order sensitivity of a response to variations in the model parameter is the G-differential* $\delta R(\mathbf{e}^0; \mathbf{h})$ defined in Eq. (6). Since the system's state vector $\mathbf{u}$ and parameters $\boldsymbol{\alpha}$ are related to each other through Eqs. (1) and (2), it follows that $\mathbf{h}_u$ and $\mathbf{h}_\alpha$ are also related to each other. Therefore, the sensitivity $\delta R(\mathbf{e}^0; \mathbf{h})$ of $R(\mathbf{e})$ at $\mathbf{e}^0$ can only be evaluated after determining the vector of variations $\mathbf{h}_u$ in terms of the vector of parameter variations $\mathbf{h}_\alpha$. The first-order relationship between $\mathbf{h}_u$ and $\mathbf{h}_\alpha$ is determined by taking the G-differentials of Eqs. (1) and (2). Thus, taking the G-differential at $\mathbf{e}^0$ of Eq. (1) yields the following operator-block-matrix equation:

$$\delta \mathbf{N}\left(\mathbf{u}^{0}, \boldsymbol{\alpha}^{0}; \mathbf{h}_{u}, \mathbf{h}_{\alpha}\right) = \delta \mathbf{Q}\left(\boldsymbol{\alpha}^{0}; \mathbf{h}_{\alpha}\right), \quad \mathbf{x} \in \Omega_{x}. \tag{7}$$

The boundary conditions associated with Eq. (2.7) are obtained by taking the G-differential at $\mathbf{e}^{0}$ of the boundary and initial conditions represented by Eq. (2.2), which yields

$$\delta \mathbf{B}\left(\mathbf{u}^{0}, \boldsymbol{\alpha}^{0}; \mathbf{h}_{u}, \mathbf{h}_{\alpha}\right) = \delta \mathbf{A}\left(\boldsymbol{\alpha}^{0}; \mathbf{h}_{\alpha}\right), \quad \mathbf{x} \in \partial\Omega_{x}. \tag{8}$$

Equations (7) and (8) represent the "*first-level forward sensitivity system*" ($1^{st}$-LFSS). For a given vector of parameter variations $\mathbf{h}_{\alpha}$ around $\boldsymbol{\alpha}^{0}$, the $1^{st}$-LFSS represented by Eqs. (7) and (8) is solved to obtain $\mathbf{h}_{u}$. Once $\mathbf{h}_{u}$ is available, it is in turn used in Eq. (6) to compute the sensitivity $\delta R(\mathbf{e}^{0}; \mathbf{h})$ of $R(\mathbf{e})$ at $\mathbf{e}^{0}$, for a given vector of parameter variations $\mathbf{h}_{\alpha}$. *The direct computation of the response sensitivity $\delta R(\mathbf{e}^{0}; \mathbf{h})$ by using the ($\mathbf{h}_{\alpha}$-dependent) solution $\mathbf{h}_{u}$ of the $1^{st}$-LFSS is called* [1] *the (first-order) forward sensitivity analysis method ($1^{st}$-FSAM).* As is well known, the $1^{st}$-FSAM requires $O(N_{\alpha})$ large-scale forward computations; therefore, the $1^{st}$-FSAM is advantageous to employ only if, in the problem under consideration, the number $N_{r}$ of responses of interest exceeds the number of system parameters and/or parameter variations of interest.

In most practical situations, however, the number of model parameters exceeds significantly the number of functional responses of interest, i.e., $N_{r} \ll N_{\alpha}$. In such cases, *the (first-order) adjoint sensitivity analysis methodology ($1^{st}$-ASAM), generally developed by Cacuci [1-2], is the most efficient method for computing exactly the first-order sensitiviti*es, *since it requires only $O(N_{r})$ large-scale computations*. For constructing the $1^{st}$-ASAM, it is necessary that $\delta R(\mathbf{e}^{0}; \mathbf{h})$ be linear in $\mathbf{h}$, which implies that $R(\mathbf{e})$ must satisfy a weak Lipschitz condition at $\mathbf{e}^{0}$, and also satisfy the following condition

$$R\left(\mathbf{e}^{0} + \varepsilon\mathbf{h}_{1} + \varepsilon\mathbf{h}_{2}\right) - R\left(\mathbf{e}^{0} + \varepsilon\mathbf{h}_{1}\right) - R\left(\mathbf{e}^{0} + \varepsilon\mathbf{h}_{2}\right) + R\left(\mathbf{e}^{0}\right) = o(t); \\ \mathbf{h}_{1}, \mathbf{h}_{2} \in \mathcal{H}_{u} \times \mathcal{H}_{\alpha}; \quad \varepsilon \in \mathcal{F}. \tag{9}$$

If $R(\mathbf{e})$ satisfies the two conditions above, then the total response variation $\delta R(\mathbf{e}^0;\mathbf{h})$ is indeed linear in $\mathbf{h}$, and can therefore be denoted as $DR(\mathbf{e}^0;\mathbf{h})$. Consequently, $R(\mathbf{e})$ admits a total G-derivative at $\mathbf{e}^0 = (\mathbf{u}^0, \boldsymbol{\alpha}^0)$, such that the relationship

$$DR(\mathbf{e}^0;\mathbf{h}) = R'_u(\mathbf{e}^0)\mathbf{h}_u + R'_\alpha(\mathbf{e}^0)\mathbf{h}_\alpha \tag{10}$$

holds, where $R'_u(\mathbf{e}^0)$ and $R'_\alpha(\mathbf{e}^0)$ denote the partial G-derivatives at $\mathbf{e}^0$ of $R(\mathbf{e})$ with respect to $\mathbf{u}$ and $\boldsymbol{\alpha}$. It is convenient to refer to the quantities $R'_u(\mathbf{e}^0)\mathbf{h}_u$ and $R'_\alpha(\mathbf{e}^0)\mathbf{h}_\alpha$ appearing in Eq. (10) as the "*indirect effect term*" and the "*direct effect term*," respectively. The operator $R'_u(\mathbf{e}^0)$ acts linearly on the vector of (arbitrary) variations $\mathbf{h}_u$, from $\mathscr{E}_u$ into $\mathscr{F}$, while the operator $R'_\alpha(\mathbf{e}^0)$ acts linearly on the vector of (arbitrary) variations $\mathbf{h}_\alpha$, from $\mathscr{E}_\alpha$ into $\mathscr{F}$.

To implement the $1^{st}$-ASAM for computing the first-order G-differential $DR(\mathbf{e}^0;\mathbf{h})$, the spaces $\mathscr{E}_u$, $\mathscr{E}_\alpha$, and $\mathscr{E}_Q$ will henceforth be considered to be Hilbert spaces and denoted as $\mathscr{H}_u(\Omega_x)$ and $\mathscr{H}_Q(\Omega_x)$, respectively. The elements of $\mathscr{H}_u(\Omega_x)$ and $\mathscr{H}_Q(\Omega_x)$ are, as before, vector-valued functions defined on the open set $\Omega_x \subset \mathbb{R}^{J_x}$, with smooth boundary $\partial\Omega_x$. On $\mathscr{H}_u(\Omega_x)$, the inner product of two vectors $\mathbf{u}^{(1)} \in \mathscr{H}_u$ and $\mathbf{u}^{(2)} \in \mathscr{H}_u$ will be denoted as $\langle \mathbf{u}^{(1)}, \mathbf{u}^{(2)} \rangle_u$; similarly, the inner product on $\mathscr{H}_Q(\Omega_x)$ of two vectors $\mathbf{Q}^{(1)} \in \mathscr{H}_Q$ and $\mathbf{Q}^{(2)} \in \mathscr{H}_Q$ will be denoted as $\langle \mathbf{Q}^{(1)}, \mathbf{Q}^{(2)} \rangle_Q$. The inner product of two vectors $\boldsymbol{\alpha}^{(1)} \in \mathscr{H}_\alpha$ and $\boldsymbol{\alpha}^{(2)} \in \mathscr{H}_\alpha$ on the Hilbert space $\mathscr{H}_\alpha$ will be denoted as $\langle \boldsymbol{\alpha}^{(1)}, \boldsymbol{\alpha}^{(2)} \rangle_\alpha$.

Note that $R'_u(\mathbf{e}^0)\mathbf{h}_u$ is a continuous linear (in $\mathbf{h}_u$) functional on $\mathscr{H}_u(\Omega_x)$, and thus lies in the dual space $\mathscr{H}_u^*(\Omega_x)$. As is well known, real Hilbert spaces are self-dual; for complex Hilbert spaces, it is conventional to identify $\mathscr{H}_u^*(\Omega_x)$ with $\mathscr{H}_u(\Omega_x)$ in the sense of correspondence under the usual isometric anti-isomorphism introduced by complex conjugation. Furthermore, the well known Riesz representation theorem ensures that there exists a unique (column) vector $\mathbf{D}_u R(e^0) \in \mathcal{H}_u$, defined as $\mathbf{D}_u R(e) \equiv [\partial R(e)/\partial u_1, \ldots, \partial R(e)/\partial u_{K_u}]$, which is customarily called the partial gradient of $R(e)$ with respect to $\mathbf{u}$, evaluated at $e^0$, such that

$$R'_u(e^0)\mathbf{h}_u = \langle \mathbf{D}_u R(e^0), \mathbf{h}_u \rangle_u, \quad \mathbf{h}_u \in \mathcal{H}_u. \tag{11}$$

Similarly, the functional $R'_\alpha(e^0)\mathbf{h}_\alpha$ is linear in $\mathbf{h}_\alpha$; consequently, the Riesz representation theorem ensures that there exists a unique vector $\mathbf{D}_\alpha R(e^0) \in \mathcal{H}_\alpha$, where $\mathbf{D}_\alpha R(e) \equiv [\partial R(e)/\partial \alpha_1, \ldots, \partial R(e)/\partial \alpha_{N_\alpha}]$, which is customarily called the partial gradient of $R(e)$ with respect to $\boldsymbol{\alpha}$, evaluated at $e^0$, such that

$$R'_\alpha(e^0)\mathbf{h}_\alpha = \langle \mathbf{D}_\alpha R(e^0), \mathbf{h}_\alpha \rangle_\alpha, \quad \mathbf{h}_\alpha \in \mathcal{H}_\alpha. \tag{12}$$

The *1st-ASAM* also requires that the operators $\mathbf{N}(\mathbf{u},\boldsymbol{\alpha})$ and $\mathbf{B}(\mathbf{u},\boldsymbol{\alpha})$ satisfy the same conditions as the operator $R(\mathbf{e})$, namely that $\mathbf{N}(\mathbf{u},\boldsymbol{\alpha})$ and $\mathbf{B}(\mathbf{u},\boldsymbol{\alpha})$ each satisfy a weak Lipschitz condition at $e^0$, and also satisfy the condition given in Eq. (9). Under these conditions, the operators $\delta \mathbf{N}(\mathbf{u}^0,\boldsymbol{\alpha}^0;\mathbf{h}_u,\mathbf{h}_\alpha)$ and $\delta \mathbf{B}(\mathbf{u}^0,\boldsymbol{\alpha}^0;\mathbf{h}_u,\mathbf{h}_\alpha)$ will each be separately linear in $\mathbf{h}_u$ and $\mathbf{h}_\alpha$, respectively. Finally, the *1st-ASAM* also requires that the operators $\mathbf{Q}(\boldsymbol{\alpha})$ and $\mathbf{A}(\boldsymbol{\alpha})$ each satisfy a weak Lipschitz condition at $\boldsymbol{\alpha}^0$, and also satisfy the condition given in Eq. (9) at $\boldsymbol{\alpha}^0$. Under these conditions, the operators $\delta \mathbf{Q}(\boldsymbol{\alpha}^0;\mathbf{h}_\alpha)$ and $\delta \mathbf{A}(\boldsymbol{\alpha}^0;\mathbf{h}_\alpha)$ will each be linear in $\mathbf{h}_\alpha$. Consequently, Eqs. (7) and (8) can be written, respectively, in the following inner-product forms

$$\mathbf{D}_u \mathbf{N}(\mathbf{u}^0,\boldsymbol{\alpha}^0)\mathbf{h}_u = \mathbf{D}_\alpha \mathbf{Q}(\boldsymbol{\alpha}^0)\mathbf{h}_\alpha - \mathbf{D}_\alpha \mathbf{N}(\mathbf{u}^0,\boldsymbol{\alpha}^0)\mathbf{h}_\alpha, \quad \mathbf{x} \in \Omega_x, \tag{13}$$

and

$$\mathbf{D}_u \mathbf{B}(\mathbf{u}^0,\boldsymbol{\alpha}^0)\mathbf{h}_u = \mathbf{D}_\alpha \mathbf{A}(\boldsymbol{\alpha}^0)\mathbf{h}_\alpha - \mathbf{D}_\alpha \mathbf{B}(\mathbf{u}^0,\boldsymbol{\alpha}^0)\mathbf{h}_\alpha, \quad \mathbf{x} \in \partial\Omega_x, \tag{14}$$

respectively, where $\mathbf{D}_u \mathbf{N}(\mathbf{u}^0,\boldsymbol{\alpha}^0)$ and $\mathbf{D}_\alpha \mathbf{N}(\mathbf{u}^0,\boldsymbol{\alpha}^0)$ denote the respective partial G-derivatives, defined as

$$\mathbf{D}_u \mathbf{N}(\mathbf{u}^0, \boldsymbol{\alpha}^0) \equiv \begin{pmatrix} \dfrac{\partial N_1(\mathbf{u}, \boldsymbol{\alpha})}{\partial u_1} & \cdots & \dfrac{\partial N_1(\mathbf{u}, \boldsymbol{\alpha})}{\partial u_{N_u}} \\ \vdots & \ddots & \vdots \\ \dfrac{\partial N_{K_u}(\mathbf{u}, \boldsymbol{\alpha})}{\partial u_1} & \cdots & \dfrac{\partial N_{K_u}(\mathbf{u}, \boldsymbol{\alpha})}{\partial u_{N_u}} \end{pmatrix}, \quad \mathbf{D}_\alpha \mathbf{N}(\mathbf{u}^0, \boldsymbol{\alpha}^0) \equiv \begin{pmatrix} \dfrac{\partial N_1(\mathbf{u}, \boldsymbol{\alpha})}{\partial \alpha_1} & \cdots & \dfrac{\partial N_1(\mathbf{u}, \boldsymbol{\alpha})}{\partial \alpha_{N_\alpha}} \\ \vdots & \ddots & \vdots \\ \dfrac{\partial N_{K_u}(\mathbf{u}, \boldsymbol{\alpha})}{\partial \alpha_1} & \cdots & \dfrac{\partial N_{K_u}(\mathbf{u}, \boldsymbol{\alpha})}{\partial \alpha_{N_\alpha}} \end{pmatrix}.$$

The partial G-derivatives $\mathbf{D}_u \mathbf{B}(\mathbf{u}^0, \boldsymbol{\alpha}^0)$, $\mathbf{D}_\alpha \mathbf{Q}(\boldsymbol{\alpha}^0)$, $\mathbf{D}_\alpha \mathbf{A}(\boldsymbol{\alpha}^0)$, and $\mathbf{D}_\alpha \mathbf{B}(\boldsymbol{\alpha}^0)\mathbf{u}^0$ are defined analogously to the ones above.

Following [1] for the construction of the $1^{st}$-ASAM, the formal adjoint $\mathbf{L}^*(\boldsymbol{\alpha}^0)$ of $\mathbf{L}(\boldsymbol{\alpha}^0)$ is now introduced by requiring that the following relationship hold for an arbitrary vector $\boldsymbol{\psi}(\mathbf{x}) \equiv [\psi_1(\mathbf{x}), \ldots, \psi_{N_u}(\mathbf{x})] \in \mathcal{H}_Q$:

$$\left\langle \boldsymbol{\psi}, \mathbf{D}_u \mathbf{N}(\mathbf{u}^0, \boldsymbol{\alpha}^0) \mathbf{h}_u \right\rangle_Q = \left\langle \mathbf{N}^{(1,F)}(\mathbf{u}^0, \boldsymbol{\alpha}^0) \boldsymbol{\psi}, \mathbf{h}_u \right\rangle_u + \{P(\mathbf{h}_u; \boldsymbol{\psi})\}_{\partial \Omega_x} \tag{15}$$

In the above equation, $\mathbf{N}^{(1,F)}(\mathbf{u}^0, \boldsymbol{\alpha}^0)$ denotes the *formal adjoint operator* to $\mathbf{D}_u \mathbf{N}(\mathbf{u}^0, \boldsymbol{\alpha}^0)$, and is the $N_u \times N_u$ matrix obtained by transposing the formal adjoints of components of $\mathbf{D}_u \mathbf{N}(\mathbf{u}^0, \boldsymbol{\alpha}^0)$. Furthermore in Eq. (15), the quantity $\{P(\mathbf{h}_u; \boldsymbol{\psi})\}_{\partial \Omega_x}$ denotes the *associated bilinear form* evaluated on $\partial \Omega_x$. The domain of $\mathbf{N}^{(1,F)}(\mathbf{u}^0, \boldsymbol{\alpha}^0)$ is determined by selecting appropriate *adjoint boundary and/or initial conditions*, which are obtained by using Eqs. (13) and (14), and by requiring that they (i.e., the *adjoint boundary and/or initial conditions*) must be independent of *unknown* values of $\mathbf{h}_u$ and $\mathbf{h}_\alpha$. In certain situations, might be computationally advantageous to include certain boundary components of $\{P(\mathbf{h}_u; \boldsymbol{\psi})\}_{\partial \Omega_x}$ into the components of $\mathbf{N}^{(1,F)}(\mathbf{u}^0, \boldsymbol{\alpha}^0)$, in which case the resulting *extended adjoint operator*, which will be denoted as $\mathbf{N}^{(1)}(\mathbf{u}^0, \boldsymbol{\alpha}^0)$, might differ somewhat from the formal adjoint, $\mathbf{N}^{(1,F)}(\mathbf{u}^0, \boldsymbol{\alpha}^0)$, of $\mathbf{D}_u \mathbf{N}(\mathbf{u}^0, \boldsymbol{\alpha}^0)$. The *adjoint boundary and/or initial conditions* thus obtained are represented in operator form as

$$\mathbf{B}^{(1)}\left(\mathbf{u}^0, \boldsymbol{\alpha}^0; \boldsymbol{\psi}\right) = \mathbf{0}, \quad \mathbf{x} \in \partial\Omega_x. \tag{16}$$

The above boundary conditions are usually inhomogeneous, i.e., $\mathbf{B}^{(1)}\left(\mathbf{u}^0, \boldsymbol{\alpha}^0; \mathbf{0}\right) \neq \mathbf{0}$. Using both the forward and the adjoint boundary and/or initial coditions represented by Eqs. (14) and (16), respectively, in Eq. (15) reduces the bilinear concomitant $\left\{P(\mathbf{h}_u; \boldsymbol{\psi})\right\}_{\partial\Omega_x}$ to a quantity, denoted as $\hat{P}\left(\mathbf{h}_\alpha, \boldsymbol{\psi}; \boldsymbol{\alpha}^0\right)$, which will contain boundary terms involving only known values of $\mathbf{h}_\alpha$, $\boldsymbol{\psi}$, and, possibly, $\boldsymbol{\alpha}^0$. In general, $\hat{P}$ does not automatically vanish as a result of the operations discussed in the foregoing, although it may do so in particular instances. The result of these operations for choosing the most advantageos (for subsequent computations) adjoint boundary conditions and using Eq. (13) will transform Eq. (15) into the form

$$\left\langle \mathbf{N}^{(1)}\left(\mathbf{u}^0, \boldsymbol{\alpha}^0\right)\boldsymbol{\psi}, \mathbf{h}_u \right\rangle_u = \left\langle \boldsymbol{\psi}, \left\langle \mathbf{D}_\alpha \mathbf{Q}\left(\boldsymbol{\alpha}^0\right) - \mathbf{D}_\alpha \mathbf{N}\left(\mathbf{u}^0, \boldsymbol{\alpha}^0\right), \mathbf{h}_\alpha \right\rangle_\alpha \right\rangle_Q - \hat{P}\left(\mathbf{h}_\alpha, \boldsymbol{\psi}; \boldsymbol{\alpha}^0\right). \tag{17}$$

Since $\boldsymbol{\psi}$ is not completely defined yet, we now complete its definition by requiring that the left-side of Eq. (17) and the right-side of Eq. (11) represent the same functional, i.e.,

$$\left\langle \mathbf{N}^{(1)}\left(\mathbf{u}^0, \boldsymbol{\alpha}^0\right)\boldsymbol{\psi}, \mathbf{h}_u \right\rangle_u = \left\langle \mathbf{D}_u R\left(e^0\right), \mathbf{h}_u \right\rangle_u, \quad \mathbf{h}_u \in \mathcal{H}_u, \tag{18}$$

which implies that the adjoint function $\boldsymbol{\psi}$ is the *weak solution* of

$$\mathbf{N}^{(1)}\left(\mathbf{u}^0, \boldsymbol{\alpha}^0\right)\boldsymbol{\psi} = \mathbf{D}_u R\left(e^0\right). \tag{19}$$

Note that the well-known Riesz representation theorem ensures that the above relationship shown in Eq. (18), where $\boldsymbol{\psi}$ satisfies the adjoint boundary conditions given in Eq. (16), holds uniquely. The construction of the *first-level adjoint sensitivity system* ($1^{st}$-*LASS*), consisting of Eqs. (19) and (16), has thus been accomplished. Furthermore, Eqs. (10) through (19) can now be used to obtain the following expression for the total sensitivity $DR(\mathbf{e}^0; \mathbf{h})$ of $R(\mathbf{e})$ at:

$$DR(\mathbf{e}^0;\mathbf{h}) = \left\langle \mathbf{D}_\alpha R(e^0), \mathbf{h}_\alpha \right\rangle_\alpha + \left\langle \psi, \left\langle \mathbf{D}_\alpha \mathbf{Q}(\boldsymbol{\alpha}^0) - \mathbf{D}_\alpha \mathbf{N}(\mathbf{u}^0, \boldsymbol{\alpha}^0), \mathbf{h}_\alpha \right\rangle_\alpha \right\rangle_Q - \hat{P}(\mathbf{h}_\alpha, \psi^0; \boldsymbol{\alpha}^0)$$

$$\equiv \mathbf{S}(\mathbf{u}^0, \psi^0, \boldsymbol{\alpha}^0) \mathbf{h}_\alpha = \sum_{i=1}^{N_\alpha} S_i(\mathbf{u}^0, \psi^0, \boldsymbol{\alpha}^0) \delta \alpha_i,$$

(20)

where $\mathbf{S}(\mathbf{u}^0, \psi^0, \boldsymbol{\alpha}^0) \equiv (S_1, ..., S_{N_\alpha})^\dagger$, and where the $i^{th}$-partial first-order sensitivity (G-derivative), $S_i(\mathbf{u}^0, \psi^0, \boldsymbol{\alpha}^0)$, of $R(\mathbf{e})$ with respect to the $i^{th}$-model parameter $\alpha_i$, $i = 1, ..., N_\alpha$, is given by the expression

$$S_i(\mathbf{u}^0, \psi^0, \boldsymbol{\alpha}^0) \equiv \frac{\partial R(\mathbf{u}^0, \boldsymbol{\alpha}^0)}{\partial \alpha_i} + \left\langle \psi^0, \frac{\partial \mathbf{Q}(\mathbf{u}^0, \boldsymbol{\alpha}^0)}{\partial \alpha_i} - \frac{\partial \mathbf{N}(\mathbf{u}^0, \boldsymbol{\alpha}^0)}{\partial \alpha_i} \right\rangle_\psi - \frac{\partial \hat{P}(\mathbf{u}^0, \psi^0, \boldsymbol{\alpha}^0)}{\partial \alpha_i}, \quad i = 1, ..., N_\alpha.$$

(21)

All partial derivatives in the above expressions are to be understood as partial G-derivatives, of course. As Eq. (21) indicates, the desired elimination of all unknown values of $\mathbf{h}_u$ from the expressions of the sensitivities $S_i(\mathbf{u}, \psi, \boldsymbol{\alpha})$, $i = 1, ..., N_\alpha$, of $R(\mathbf{e})$ at $\mathbf{e}^0$ has been accomplished. Note that we have designated the space $\mathcal{H}_Q(\Omega_x)$ as $\mathcal{H}_\psi(\Omega_x)$, in order to emphasize that we are dealing with the Hilbert space on which the adjoint function $\psi$ is defined. *The sensitivities $S_i(\mathbf{u}, \boldsymbol{\alpha}, \psi)$ can therefore be computed by means of Eq. (21), after solving **only once** the $1^{st}$-level adjoint sensitivity system ($1^{st}$-LASS), consisting of Eqs. (16) and (19), to obtain the adjoint function $\psi$*. It is very important to note the $1^{st}$-LASS is independent of the functions $\mathbf{h}_u$.

Furthermore, it can be shown by straightforward computations that, *when Eqs. (1) and (2) are linear in $\mathbf{u}(\mathbf{x})$, i.e., when $\mathbf{N}[\mathbf{u}(\mathbf{x}), \boldsymbol{\alpha}(\mathbf{x})] \equiv \mathbf{L}[\boldsymbol{\alpha}(\mathbf{x})]\mathbf{u}(\mathbf{x})$ and $\mathbf{B}[\mathbf{u}(\mathbf{x}), \boldsymbol{\alpha}(\mathbf{x})] \equiv \mathbf{B}[\boldsymbol{\alpha}(\mathbf{x})]\mathbf{u}(\mathbf{x})$, then the $1^{st}$-LASS, i.e., Eqs. (16) and (19)* will reduce to the corresponding equations for linear systems [see Eqs. (18b) and (16), respectively, in Cacuci (2015)]. *The $1^{st}$-LASS for linear systems is thus independent not only of the functions $\mathbf{h}_u$ but also of the nominal values $\mathbf{u}^0$ of $\mathbf{u}$*. As discussed by Cacuci (2015), *the $1^{st}$-LASS for linear* *systems can be solved independently of the solution* $\mathbf{u}^0$ *of the original equations*, which simplifies considerably the choice of numerical methods for solving the adjoint system.

## 3. THE SECOND-ORDER ADJOINT SENSITIVITY ANALYSIS METHODOLOGY (2$^{nd}$-ASAM) FOR NONLINEAR SYSTEMS

We now note the very important fact that, since Eq. (21) holds for any nominal parameter values, $\boldsymbol{\alpha}^0$, it follows that the first-order sensitivities $S_i(\mathbf{u}, \boldsymbol{\psi}, \boldsymbol{\alpha})$, $i = 1, \ldots, N_\alpha$, can be generally considered as functionals of the original state-function $\mathbf{u}$, the parameters $\boldsymbol{\alpha}$, and the adjoint function $\boldsymbol{\psi}$, namely

$$S_i(\mathbf{u}, \boldsymbol{\psi}, \boldsymbol{\alpha}) = \frac{\partial R(\mathbf{u}, \boldsymbol{\alpha})}{\partial \alpha_i} + \left\langle \boldsymbol{\psi}, \frac{\partial \mathbf{Q}(\boldsymbol{\alpha})}{\partial \alpha_i} - \frac{\partial \mathbf{N}(\mathbf{u}, \boldsymbol{\alpha})}{\partial \alpha_i} \right\rangle_\psi - \frac{\partial \hat{P}(\mathbf{u}, \boldsymbol{\psi}, \boldsymbol{\alpha})}{\partial \alpha_i}, \quad i = 1, \ldots, N_\alpha, \quad (22)$$

where the original state-function $\mathbf{u}$ satisfies Eqs. (1) and (2) while the adjoint function $\boldsymbol{\psi}$ satisfies the *1$^{st}$-LASS*

$$\mathbf{N}^{(1)}(\mathbf{u}, \boldsymbol{\alpha}) \boldsymbol{\psi} = \mathbf{D}_u R(\mathbf{u}, \boldsymbol{\alpha}), \quad \mathbf{x} \in \Omega_x, \quad (23)$$

$$\mathbf{B}^{(1)}(\mathbf{u}, \boldsymbol{\alpha}; \boldsymbol{\psi}) = \mathbf{0}, \quad \mathbf{x} \in \partial \Omega_x., \quad \mathbf{x} \in \partial \Omega_x. \quad (24)$$

As Eq. (22) indicates, the first-order sensitivities are functionals of the form $S_i(\mathbf{u}, \boldsymbol{\psi}, \boldsymbol{\alpha}) : \mathcal{D}_{R_i} \subset \mathcal{H}_u(\Omega_x) \times \mathcal{H}_\psi(\Omega_x) \times \mathcal{H}_\alpha \to \mathcal{F}$. Hence, it is possible to define the first-order G-differential, $\delta S_i(\mathbf{u}^0, \boldsymbol{\psi}^0, \boldsymbol{\alpha}^0; \mathbf{h}_u, \mathbf{h}_\psi, \mathbf{h}_\alpha)$, of any of the functionals $S_i(\mathbf{u}, \boldsymbol{\psi}, \boldsymbol{\alpha})$, at an arbitrary point $(\mathbf{e}^0, \boldsymbol{\psi}^0)$, in the usual manner, namely

$$\delta S_i(\mathbf{u}^0, \boldsymbol{\psi}^0, \boldsymbol{\alpha}^0; \mathbf{h}_u, \mathbf{h}_\psi, \mathbf{h}_\alpha) \equiv \left\{ \frac{d}{d\varepsilon} \left[ S_i(\mathbf{u}^0 + \varepsilon \mathbf{h}_u, \boldsymbol{\psi}^0 + \varepsilon \mathbf{h}_\psi, \boldsymbol{\alpha}^0 + \varepsilon \mathbf{h}_\alpha) \right] \right\}_{\varepsilon=0}, \quad \varepsilon \in \mathcal{F}, \quad (25)$$

for all (i.e., arbitrary) vectors $(\mathbf{h}_u, \mathbf{h}_\alpha, \mathbf{h}_\psi) \in \mathcal{H}_u(\Omega_x) \times \mathcal{H}_\alpha \times \mathcal{H}_\psi(\Omega_x)$. Applying the above definition to the expression of $S_i(\mathbf{u}, \boldsymbol{\alpha}, \boldsymbol{\psi})$ given by Eq. (22) yields

$$\delta S_i\left(\mathbf{u}^0, \mathbf{\psi}^0, \mathbf{\alpha}^0; \mathbf{h}_u, \mathbf{h}_\psi, \mathbf{h}_\alpha\right) = \left(\delta S_i\right)_{direct} + \left(\delta S_i\right)_{indirect}, \tag{26}$$

where $\left(\delta S_i\right)_{direct}$ denotes the "*direct-effect term*" and $\left(\delta S_i\right)_{indirect}$ denotes the "*indirect-effect term.*" The "*direct-effect term,*" $\left(\delta S_i\right)_{direct}$, is defined as

$$\begin{aligned}\left(\delta S_i\right)_{direct} &\equiv \left\langle \mathbf{D}_\alpha S_i, \mathbf{h}_\alpha \right\rangle_\alpha \\ &= \frac{\partial}{\partial \alpha_i}\left\{\left\langle \mathbf{D}_\alpha R(\mathbf{u}, \mathbf{\alpha}) - \mathbf{D}_\alpha \hat{P}(\mathbf{u}, \mathbf{\psi}, \mathbf{\alpha}) + \left\langle \mathbf{\psi}, \mathbf{D}_\alpha Q(\mathbf{\alpha}) - \mathbf{D}_\alpha N(\mathbf{u}, \mathbf{\alpha}) \right\rangle_\psi, \mathbf{h}_\alpha \right\rangle_\alpha \right\}_{(\mathbf{u}^0, \mathbf{\psi}^0; \mathbf{\alpha}^0)},\end{aligned} \tag{27}$$

and *can be computed immediately* at this stage, without needing any additional large-scale computations. The "*indirect-effect term*", $\left(\delta S_i\right)_{indirect}$, is defined as

$$\left(\delta S_i\right)_{indirect} \equiv \left\langle \left(\mathbf{h}_u^\dagger, \mathbf{h}_\psi^\dagger\right) \begin{pmatrix} \mathbf{D}_u S_i\left(\mathbf{u}^0, \mathbf{\psi}^0, \mathbf{\alpha}^0\right) \\ \mathbf{D}_\psi S_i\left(\mathbf{u}^0, \mathbf{\psi}^0, \mathbf{\alpha}^0\right) \end{pmatrix} \right\rangle_{u\psi}, \tag{28}$$

where

$$\mathbf{D}_\psi S_i(\mathbf{u}, \mathbf{\psi}, \mathbf{\alpha}) \equiv \frac{\partial}{\partial \alpha_i}\left[Q(\mathbf{\alpha}) - N(\mathbf{u}, \mathbf{\alpha}) - \mathbf{D}_\psi \hat{P}(\mathbf{u}, \mathbf{\psi}, \mathbf{\alpha})\right], \tag{29}$$

$$\mathbf{D}_u S_i(\mathbf{u}, \mathbf{\psi}, \mathbf{\alpha}) \equiv \frac{\partial}{\partial \alpha_i}\left[\mathbf{D}_u R(\mathbf{u}, \mathbf{\alpha}) - \mathbf{D}_u \hat{P}(\mathbf{u}, \mathbf{\psi}, \mathbf{\alpha}) - \left\langle \mathbf{\psi}, \mathbf{D}_u N(\mathbf{u}, \mathbf{\alpha}) \right\rangle_\psi\right]. \tag{30}$$

Note that the "indirect-effect term" *cannot be computed at this stage*, since the vectors of variations $\mathbf{h}_u$ and $\mathbf{h}_\psi$ are unknown. Recall that the vector $\mathbf{h}_u$ is the $\mathbf{h}_\alpha$–dependent solution of the *1$^{st}$-LFSS*, comprising Eqs. (7) and (8), which are computationally impractical to solve for large-scale systems with many parameters. Furthermore, the vector of variations $\mathbf{h}_\psi$ (around the nominal value $\mathbf{\psi}^0$) is the solution of the system of equations that results from applying the definition of the G-differential to the *1$^{st}$-LASS*, i.e., to Eqs. (23) and (24), to obtain

$$\begin{aligned}&\left[\mathbf{D}_u \mathbf{N}^{(1)}\left(\mathbf{u}^0, \mathbf{\alpha}^0\right)\mathbf{\psi}^0\right]\mathbf{h}_u + \mathbf{N}^{(1)}\left(\mathbf{u}^0, \mathbf{\alpha}^0\right)\mathbf{h}_\psi - \mathbf{D}_{u^2}^2 R\left(\mathbf{u}^0, \mathbf{\alpha}^0\right)\mathbf{h}_u \\ &= \mathbf{D}_{u\alpha}^2 R\left(\mathbf{u}^0, \mathbf{\alpha}^0\right)\mathbf{h}_\alpha - \left[\mathbf{D}_\alpha \mathbf{N}^{(1)}\left(\mathbf{u}^0, \mathbf{\alpha}^0\right)\mathbf{\psi}^0\right]\mathbf{h}_\alpha, \quad \mathbf{x} \in \Omega_x,\end{aligned} \tag{31}$$

$$\left[\mathbf{D}_u\mathbf{B}^{(1)}\left(\mathbf{u}^0,\boldsymbol{\alpha}^0\right)\boldsymbol{\psi}^0\right]\mathbf{h}_u + \mathbf{B}^{(1)}\left(\mathbf{u}^0,\boldsymbol{\alpha}^0\right)\mathbf{h}_\psi = -\left[\mathbf{D}_\alpha\mathbf{B}^{(1)}\left(\mathbf{u}^0,\boldsymbol{\alpha}^0\right)\boldsymbol{\psi}^0\right]\mathbf{h}_\alpha, \quad \mathbf{x}\in\partial\Omega_x. \quad (32)$$

The operators $\mathbf{D}_{u^2}^2 R(e^0)$ and $\mathbf{D}_{\alpha u}^2 R(e^0)$, which appear in Eq. (31), are matrices of partial G-derivatives of the form

$$\mathbf{D}_{\alpha u}^2 R(e) \equiv \begin{pmatrix} \dfrac{\partial^2 R(e)}{\partial u_1 \partial \alpha_1} & \cdots & \dfrac{\partial^2 R(e)}{\partial u_1 \partial \alpha_{N_\alpha}} \\ \vdots & \ddots & \vdots \\ \dfrac{\partial^2 R(e)}{\partial u_{N_u} \partial \alpha_1} & \cdots & \dfrac{\partial^2 R(e)}{\partial u_{N_u} \partial \alpha_{N_\alpha}} \end{pmatrix}; \quad \mathbf{D}_{u^2}^2 R(e) \equiv \begin{pmatrix} \dfrac{\partial^2 R(e)}{\partial u_1^2} & \cdots & \dfrac{\partial^2 R(e)}{\partial u_1 \partial u_{N_u}} \\ \vdots & \ddots & \vdots \\ \dfrac{\partial^2 R(e)}{\partial u_{N_u} \partial u_1} & \cdots & \dfrac{\partial^2 R(e)}{\partial u_{N_u}^2} \end{pmatrix}. \quad (33)$$

As indicated by Eqs. (31) and (32), the vector of variations $\mathbf{h}_\psi$ (around the nominal value $\boldsymbol{\psi}^0$) is related to the vector of parameters variations $\mathbf{h}_\alpha$. Together, the *1$^{st}$-LFSS*, namely Eqs. (13) and (14), and the G-differentiated *1$^{st}$-LASS*, namely Eqs. (31) and (32), can be written in the following block-matrix-operator form:

$$\begin{pmatrix} \mathbf{D}_u \mathbf{N}\left(\mathbf{u}^0,\boldsymbol{\alpha}^0\right) & \mathbf{0} \\ \mathbf{D}_u \mathbf{N}^{(1)}\left(\mathbf{u}^0,\boldsymbol{\alpha}^0\right)\boldsymbol{\psi}^0 - \mathbf{D}_{u^2}^2 R\left(\mathbf{u}^0,\boldsymbol{\alpha}^0\right) & \mathbf{N}^{(1)}\left(\mathbf{u}^0,\boldsymbol{\alpha}^0\right) \end{pmatrix} \begin{pmatrix} \mathbf{h}_u \\ \mathbf{h}_\psi \end{pmatrix}$$
$$= \begin{pmatrix} \mathbf{D}_\alpha \mathbf{Q}\left(\boldsymbol{\alpha}^0\right)\mathbf{h}_\alpha - \mathbf{D}_\alpha \mathbf{N}\left(\mathbf{u}^0,\boldsymbol{\alpha}^0\right)\mathbf{h}_\alpha \\ \mathbf{D}_{u\alpha}^2 R\left(\mathbf{u}^0,\boldsymbol{\alpha}^0\right)\mathbf{h}_\alpha - \left[\mathbf{D}_\alpha \mathbf{N}^{(1)}\left(\mathbf{u}^0,\boldsymbol{\alpha}^0\right)\boldsymbol{\psi}^0\right]\mathbf{h}_\alpha \end{pmatrix}, \quad \mathbf{x}\in\Omega_x, \quad (34)$$

together with the corresponding G-differentiated boundary and/or initial conditions

$$\begin{pmatrix} \mathbf{D}_u \mathbf{B}\left(\mathbf{u}^0,\boldsymbol{\alpha}^0\right) & \mathbf{0} \\ \mathbf{D}_u \mathbf{B}^{(1)}\left(\mathbf{u}^0,\boldsymbol{\alpha}^0\right)\boldsymbol{\psi}^0 & \mathbf{B}^{(1)}\left(\mathbf{u}^0,\boldsymbol{\alpha}^0\right) \end{pmatrix} \begin{pmatrix} \mathbf{h}_u \\ \mathbf{h}_\psi \end{pmatrix} = \begin{pmatrix} \mathbf{D}_\alpha \mathbf{A}\left(\boldsymbol{\alpha}^0\right)\mathbf{h}_\alpha - \mathbf{D}_\alpha \mathbf{B}\left(\mathbf{u}^0,\boldsymbol{\alpha}^0\right)\mathbf{h}_\alpha \\ -\left[\mathbf{D}_\alpha \mathbf{B}^{(1)}\left(\mathbf{u}^0,\boldsymbol{\alpha}^0\right)\boldsymbol{\psi}^0\right]\mathbf{h}_\alpha \end{pmatrix}, \quad \mathbf{x}\in\partial\Omega_x. \quad (35)$$

The block matrix Eq. (34) together with the boundary and/or initial conditions represented by Eq. (35) constitute the *second-level forward sensitivity system (2$^{nd}$-LFSS)*. These equations can be solved to obtain the vectors $\mathbf{h}_u$ and $\mathbf{h}_\psi$, which can, in turn, be used in Eq. (27) to compute the indirect-effect term, $(\delta S_i)_{indirect}$. As shown in Eq. (26), this indirect effect term would then be added together with the already computed direct-effect term, $(\delta S_i)_{direct}$, to obtain the vector, $\delta S_i(\mathbf{u}^0, \boldsymbol{\psi}^0, \boldsymbol{\alpha}^0; \mathbf{h}_u, \mathbf{h}_\psi, \mathbf{h}_\alpha)$, of partial mixed second-order sensitivities of the form $\{\partial^2 R(\mathbf{u}, \boldsymbol{\alpha})/\partial \alpha_j \partial \alpha_i\}_{(\mathbf{u}^0, \boldsymbol{\alpha}^0)}$, $j = 1, ..., N_\alpha$. Thus, Eqs. (34) and (35) would be solved to obtain $\delta S_1(\mathbf{u}^0, \boldsymbol{\psi}^0, \boldsymbol{\alpha}^0; \mathbf{h}_u, \mathbf{h}_\psi, \mathbf{h}_\alpha)$, *for* $i = 1$; the process of solving these operator equations would continue for $\delta S_2(\mathbf{u}^0, \boldsymbol{\psi}^0, \boldsymbol{\alpha}^0; \mathbf{h}_u, \mathbf{h}_\psi, \mathbf{h}_\alpha)$, *for* $i = 2$, and would end by computing $\delta S_{N_\alpha}(\mathbf{u}^0, \boldsymbol{\psi}^0, \boldsymbol{\alpha}^0; \mathbf{h}_u, \mathbf{h}_\psi, \mathbf{h}_\alpha)$, *for* $i = N_\alpha$.

The computational process just described would yield the complete set of all second-order responses sensitivities, $\{\partial^2 R(\mathbf{u}, \boldsymbol{\alpha})/\partial \alpha_j \partial \alpha_i\}_{(\mathbf{u}^0, \boldsymbol{\alpha}^0)}$, to all of the system parameters, in $O(N_\alpha^2)$ large-scale forward computations. *However, such a computational burden would be impractical for large-scale systems.*

Just as the *1$^{st}$-ASAM* aimed at computing efficiently (using adjoint functions) the first-order indirect-effect term $R'_u(\mathbf{e}^0)\mathbf{h}_u$, the *second-order adjoint sensitivity analysis methodology* (*2$^{nd}$-ASAM*), to be presented in the following, aims at the efficient computation of the second-order indirect-effect term $(\delta S_i)_{indirect}$, by using the *2$^{nd}$-level adjoint sensitivity system* (*2$^{nd}$-LASS*) that will be developed next. To construct the *2$^{nd}$-LASS*, consider the Hilbert space $\mathcal{H}_{u\psi}(\Omega_x) \equiv \mathcal{H}_u(\Omega_x) \times \mathcal{H}_\psi(\Omega_x)$, comprising elements (vectors) of the form $\boldsymbol{\varphi}^{(2)} \equiv (\boldsymbol{\varphi}_1^{(2)}, \boldsymbol{\varphi}_2^{(2)})$ and $\boldsymbol{\psi}^{(2)} \equiv (\boldsymbol{\psi}_1^{(2)}, \boldsymbol{\psi}_2^{(2)})$, respectively, and endowed with the inner product $\langle \boldsymbol{\varphi}^{(2)}, \boldsymbol{\psi}^{(2)} \rangle_{u\psi}$ defined as

$$\langle \boldsymbol{\varphi}^{(2)}, \boldsymbol{\psi}^{(2)} \rangle_{u\psi} \equiv \int_{\Omega_x} \left( \boldsymbol{\varphi}_1^{(2)} \bullet \boldsymbol{\psi}_1^{(2)} + \boldsymbol{\varphi}_2^{(2)} \bullet \boldsymbol{\psi}_2^{(2)} \right) d\mathbf{x}. \tag{36}$$

Using the above definition, we form the inner product of Eq. (34) with a yet undefined vector $\boldsymbol{\psi}_i^{(2)} \equiv (\boldsymbol{\psi}_{i1}^{(2)}, \boldsymbol{\psi}_{i2}^{(2)})$ to obtain the sequence of equalities shown below:

$$\left\langle \left( \mathbf{\psi}_{i1}^{(2)\dagger}, \mathbf{\psi}_{i2}^{(2)\dagger} \right) \begin{pmatrix} \mathbf{D}_u \mathbf{N}\left(\mathbf{u}^0, \mathbf{\alpha}^0\right) & \mathbf{0} \\ \mathbf{D}_u \mathbf{N}^{(1)}\left(\mathbf{u}^0, \mathbf{\alpha}^0\right) \mathbf{\psi}^0 - \mathbf{D}_{u^2}^2 R\left(\mathbf{u}^0, \mathbf{\alpha}^0\right) & \mathbf{N}^{(1)}\left(\mathbf{u}^0, \mathbf{\alpha}^0\right) \end{pmatrix} \begin{pmatrix} \mathbf{h}_u \\ \mathbf{h}_\psi \end{pmatrix} \right\rangle_{u\psi}$$

$$= \left\langle \left( \mathbf{\psi}_{i1}^{(2)\dagger}, \mathbf{\psi}_{i2}^{(2)\dagger} \right) \begin{pmatrix} \mathbf{D}_\alpha \mathbf{Q}\left(\mathbf{\alpha}^0\right) \mathbf{h}_\alpha - \mathbf{D}_\alpha \mathbf{N}\left(\mathbf{u}^0, \mathbf{\alpha}^0\right) \mathbf{h}_\alpha \\ \mathbf{D}_{u\alpha}^2 R\left(\mathbf{u}^0, \mathbf{\alpha}^0\right) \mathbf{h}_\alpha - \left[ \mathbf{D}_\alpha \mathbf{N}^{(1)}\left(\mathbf{u}^0, \mathbf{\alpha}^0\right) \mathbf{\psi}^0 \right] \mathbf{h}_\alpha \end{pmatrix} \right\rangle_{u\psi} \quad (37)$$

$$= \left\langle \left( \mathbf{h}_u^\dagger, \mathbf{h}_\psi^\dagger \right) \begin{pmatrix} \mathbf{A}_{11}^{(F,2)}\left(\mathbf{u}^0, \mathbf{\alpha}^0\right) & \mathbf{A}_{12}^{(F,2)}\left(\mathbf{u}^0, \mathbf{\alpha}^0\right) \\ \mathbf{0} & \mathbf{A}_{22}^{(F,2)}\left(\mathbf{u}^0, \mathbf{\alpha}^0\right) \end{pmatrix} \begin{pmatrix} \mathbf{\psi}_{i1}^{(2)} \\ \mathbf{\psi}_{i2}^{(2)} \end{pmatrix} \right\rangle_{u\psi}$$

$$+ \left\{ P_2\left(\mathbf{u}^0, \mathbf{\alpha}^0; \mathbf{\psi}_{i1}^{(2)}, \mathbf{\psi}_{i2}^{(2)}; \mathbf{h}_u, \mathbf{h}_\psi\right) \right\}_{\partial\Omega_x}, \quad \mathbf{x} \in \Omega_x,$$

where $\mathbf{A}_{11}^{(F,2)}\left(\mathbf{u}^0, \mathbf{\alpha}^0\right)$ denotes the formal adjoint of the operator $\mathbf{D}_u \mathbf{N}\left(\mathbf{u}^0, \mathbf{\alpha}^0\right)$, $\mathbf{A}_{12}^{(F,2)}\left(\mathbf{u}^0, \mathbf{\alpha}^0\right)$ denotes the formal adjoint of the operator $\mathbf{D}_u \mathbf{N}^{(1)}\left(\mathbf{u}^0, \mathbf{\alpha}^0\right) \mathbf{\psi}^0 - \mathbf{D}_{u^2}^2 R\left(\mathbf{u}^0, \mathbf{\alpha}^0\right)$, $\mathbf{A}_{22}^{(F,2)}\left(\mathbf{u}^0, \mathbf{\alpha}^0\right)$ denotes the formal adjoint of the operator $\mathbf{N}^{(1)}\left(\mathbf{u}^0, \mathbf{\alpha}^0\right)$, and $\left\{ P_2\left(\mathbf{u}^0, \mathbf{\alpha}^0; \mathbf{\psi}_{i1}^{(2)}, \mathbf{\psi}_{i2}^{(2)}; \mathbf{h}_u, \mathbf{h}_\psi\right) \right\}_{\partial\Omega_x}$ denotes the corresponding bilinear concomitant on $\mathbf{x} \in \partial\Omega_x$. The respective domains of the formal adjoint operators $\mathbf{A}_{11}^{(F,2)}\left(\mathbf{u}^0, \mathbf{\alpha}^0\right)$, $\mathbf{A}_{12}^{(F,2)}\left(\mathbf{u}^0, \mathbf{\alpha}^0\right)$, and $\mathbf{A}_{22}^{(F,2)}\left(\mathbf{u}^0, \mathbf{\alpha}^0\right)$ are determined by selecting appropriate *adjoint boundary and/or initial conditions*, which are obtained by using Eq. (35), and by requiring that they (i.e., the *adjoint boundary and/or initial conditions*) must be independent of *unknown* values of $\mathbf{h}_u$, $\mathbf{h}_\psi$, and $\mathbf{h}_\alpha$. In certain situations, it might be computationally advantageous to include certain boundary components of $\left\{ P_2\left(\mathbf{u}^0, \mathbf{\alpha}^0; \mathbf{\psi}_{i1}^{(2)}, \mathbf{\psi}_{i2}^{(2)}; \mathbf{h}_u, \mathbf{h}_\psi\right) \right\}_{\partial\Omega_x}$ into the components of $\mathbf{A}_{11}^{(F,2)}\left(\mathbf{u}^0, \mathbf{\alpha}^0\right)$, $\mathbf{A}_{12}^{(F,2)}\left(\mathbf{u}^0, \mathbf{\alpha}^0\right)$, and/or $\mathbf{A}_{22}^{(F,2)}\left(\mathbf{u}^0, \mathbf{\alpha}^0\right)$. In such cases, the resulting *extended adjoint operators*, which will be denoted as $\mathbf{A}_{11}^{(2)}\left(\mathbf{u}^0, \mathbf{\alpha}^0\right)$, $\mathbf{A}_{12}^{(2)}\left(\mathbf{u}^0, \mathbf{\alpha}^0\right)$, and/or $\mathbf{A}_{22}^{(2)}\left(\mathbf{u}^0, \mathbf{\alpha}^0\right)$, might differ somewhat from the formal adjoint operators $\mathbf{A}_{11}^{(F,2)}\left(\mathbf{u}^0, \mathbf{\alpha}^0\right)$, $\mathbf{A}_{12}^{(F,2)}\left(\mathbf{u}^0, \mathbf{\alpha}^0\right)$, and/or $\mathbf{A}_{22}^{(F,2)}\left(\mathbf{u}^0, \mathbf{\alpha}^0\right)$. The *adjoint boundary and/or initial conditions* thus obtained are represented in operator form as

$$\begin{aligned} \mathbf{B}_1^{(2)}\left(\mathbf{u}^0, \mathbf{\alpha}^0; \mathbf{\psi}_{i1}^{(2)}, \mathbf{\psi}_{i2}^{(2)}\right) &= \mathbf{0}, \quad \mathbf{x} \in \partial\Omega_x, \\ \mathbf{B}_2^{(2)}\left(\mathbf{u}^0, \mathbf{\alpha}^0; \mathbf{\psi}_{i1}^{(2)}, \mathbf{\psi}_{i2}^{(2)}\right) &= \mathbf{0}, \quad \mathbf{x} \in \partial\Omega_x. \end{aligned} \quad (38)$$

The above boundary conditions are usually inhomogeneous in the functions $\boldsymbol{\psi}_i^{(2)} \equiv \left(\boldsymbol{\psi}_{i1}^{(2)}, \boldsymbol{\psi}_{i2}^{(2)}\right)$, i.e., $\mathbf{B}_1^{(2)}\left(\mathbf{u}^0, \boldsymbol{\alpha}^0; \mathbf{0}, \mathbf{0}\right) \neq \mathbf{0}$ and $\mathbf{B}_2^{(2)}\left(\mathbf{u}^0, \boldsymbol{\alpha}^0; \mathbf{0}, \mathbf{0};\right) \neq \mathbf{0}$. Using both the forward and the adjoint boundary and/or initial conditions represented by Eqs. (35) and (38), respectively, in Eq. (37) eliminates the terms containing unknown values of the vectors $\mathbf{h}_u$ and $\mathbf{h}_\psi$, reducing the bilinear concomitant $\left\{P_2\left(\mathbf{u}^0, \boldsymbol{\alpha}^0; \boldsymbol{\psi}_{i1}^{(2)}, \boldsymbol{\psi}_{i2}^{(2)}; \mathbf{h}_u, \mathbf{h}_\psi\right)\right\}_{\partial\Omega_x}$ to a quantity, denoted as $\left\{\hat{P}_2\left(\mathbf{u}^0, \boldsymbol{\alpha}^0; \boldsymbol{\psi}_{i1}^{(2)}, \boldsymbol{\psi}_{i2}^{(2)}; \mathbf{h}_\alpha\right)\right\}_{\partial\Omega_x}$, which will contain boundary terms involving only known values of $\mathbf{h}_\alpha$ and of its other arguments. In general, $\left\{\hat{P}_2\left(\mathbf{u}^0, \boldsymbol{\alpha}^0; \boldsymbol{\psi}_{i1}^{(2)}, \boldsymbol{\psi}_{i2}^{(2)}; \mathbf{h}_\alpha\right)\right\}_{\partial\Omega_x}$ does not automatically vanish as a result of the operations discussed in the foregoing, although it may do so in particular instances.

The result of the above-mentioned operations for choosing the most advantageos (for subsequent computations) adjoint boundary conditions transforms Eq. (37) into the form

$$\left\langle \left(\boldsymbol{\psi}_{i1}^{(2)\dagger}, \boldsymbol{\psi}_{i2}^{(2)\dagger}\right) \begin{pmatrix} \mathbf{D}_u \mathbf{N}\left(\mathbf{u}^0, \boldsymbol{\alpha}^0\right) & \mathbf{0} \\ \mathbf{D}_u \mathbf{N}^{(1)}\left(\mathbf{u}^0, \boldsymbol{\alpha}^0\right)\boldsymbol{\psi}^0 - \mathbf{D}_{u^2}^2 R\left(\mathbf{u}^0, \boldsymbol{\alpha}^0\right) & \mathbf{N}^{(1)}\left(\mathbf{u}^0, \boldsymbol{\alpha}^0\right) \end{pmatrix} \begin{pmatrix} \mathbf{h}_u \\ \mathbf{h}_\psi \end{pmatrix} \right\rangle_{u\psi}$$

$$= \left\langle \left(\boldsymbol{\psi}_{i1}^{(2)\dagger}, \boldsymbol{\psi}_{i2}^{(2)\dagger}\right) \begin{pmatrix} \mathbf{D}_\alpha \mathbf{Q}\left(\boldsymbol{\alpha}^0\right)\mathbf{h}_\alpha - \mathbf{D}_\alpha \mathbf{N}\left(\mathbf{u}^0, \boldsymbol{\alpha}^0\right)\mathbf{h}_\alpha \\ \mathbf{D}_{u\alpha}^2 R\left(\mathbf{u}^0, \boldsymbol{\alpha}^0\right)\mathbf{h}_\alpha - \left[\mathbf{D}_\alpha \mathbf{N}^{(1)}\left(\mathbf{u}^0, \boldsymbol{\alpha}^0\right)\boldsymbol{\psi}^0\right]\mathbf{h}_\alpha \end{pmatrix} \right\rangle_{u\psi} \tag{39}$$

$$= \left\langle \left(\mathbf{h}_u^\dagger, \mathbf{h}_\psi^\dagger\right) \begin{pmatrix} \mathbf{A}_{11}^{(2)}\left(\mathbf{u}^0, \boldsymbol{\alpha}^0\right) & \mathbf{A}_{12}^{(2)}\left(\mathbf{u}^0, \boldsymbol{\alpha}^0\right) \\ \mathbf{0} & \mathbf{A}_{22}^{(2)}\left(\mathbf{u}^0, \boldsymbol{\alpha}^0\right) \end{pmatrix} \begin{pmatrix} \boldsymbol{\psi}_{i1}^{(2)} \\ \boldsymbol{\psi}_{i2}^{(2)} \end{pmatrix} \right\rangle_{u\psi}$$

$$+ \left\{\hat{P}_2\left(\mathbf{u}^0, \boldsymbol{\alpha}^0; \boldsymbol{\psi}_1^{(2)}, \boldsymbol{\psi}_2^{(2)}; \mathbf{h}_\alpha\right)\right\}_{\partial\Omega_x}, \quad \mathbf{x} \in \Omega_x,$$

Folowing the same principles as in Section 2, we now require the first term on the right-side of the last equality in Eq. (37) to represent the same functional as the right-side of Eq. (28), i.e.,

$$\left\langle \left(\mathbf{h}_u^\dagger, \mathbf{h}_\psi^\dagger\right) \begin{pmatrix} \mathbf{A}_{11}^{(2)}\left(\mathbf{u}^0, \boldsymbol{\alpha}^0\right) & \mathbf{A}_{12}^{(2)}\left(\mathbf{u}^0, \boldsymbol{\alpha}^0\right) \\ \mathbf{0} & \mathbf{A}_{22}^{(2)}\left(\mathbf{u}^0, \boldsymbol{\alpha}^0\right) \end{pmatrix} \begin{pmatrix} \boldsymbol{\psi}_{i1}^{(2)} \\ \boldsymbol{\psi}_{i2}^{(2)} \end{pmatrix} \right\rangle_{u\psi} = \left\langle \left(\mathbf{h}_u^\dagger, \mathbf{h}_\psi^\dagger\right) \begin{pmatrix} \mathbf{D}_u S_i\left(\mathbf{u}^0, \boldsymbol{\psi}^0, \boldsymbol{\alpha}^0\right) \\ \mathbf{D}_\psi S_i\left(\mathbf{u}^0, \boldsymbol{\psi}^0, \boldsymbol{\alpha}^0\right) \end{pmatrix} \right\rangle_{u\psi} \tag{40}$$

The above relations indicate that the adjoint functions $\psi_{i1}^{(2)}(\mathbf{x})$ and $\psi_{i2}^{(2)}(\mathbf{x})$ are *weak solutions* of the system of operator equations below:

$$\begin{aligned}
\mathbf{A}_{22}^{(2)}(\mathbf{u}^0, \boldsymbol{\alpha}^0) \boldsymbol{\psi}_{i2}^{(2)} &= \mathbf{D}_\psi S_i(\mathbf{u}^0, \boldsymbol{\psi}^0, \boldsymbol{\alpha}^0), \\
\mathbf{A}_{11}^{(2)}(\mathbf{u}^0, \boldsymbol{\alpha}^0) \boldsymbol{\psi}_{i1}^{(2)} &= \mathbf{D}_u S_i(\mathbf{u}^0, \boldsymbol{\psi}^0, \boldsymbol{\alpha}^0) - \mathbf{A}_{12}^{(2)}(\mathbf{u}^0, \boldsymbol{\alpha}^0) \boldsymbol{\psi}_{i2}^{(2)},
\end{aligned} \quad (41)$$

and subject to the adjoint boundary represented by Eq. (38). After determining the adjoint functions $\psi_{i1}^{(2)}(\mathbf{x})$ and $\psi_{i2}^{(2)}(\mathbf{x})$ by solving Eqs. (41) and (38), they can be used in conjunction with Eq. (37) to represent the "*indirect-effect term*," $(\delta S_i)_{indirect}$, defined in Eq. (28) in the form

$$(\delta S_i)_{indirect} = \left\langle \left( \boldsymbol{\psi}_{i1}^{(2)\dagger}, \boldsymbol{\psi}_{i2}^{(2)\dagger} \right) \begin{pmatrix} \mathbf{D}_\alpha \mathbf{Q}(\boldsymbol{\alpha}^0) \mathbf{h}_\alpha - \mathbf{D}_\alpha \mathbf{N}(\mathbf{u}^0, \boldsymbol{\alpha}^0) \mathbf{h}_\alpha \\ \mathbf{D}_{u\alpha}^2 R(\mathbf{u}^0, \boldsymbol{\alpha}^0) \mathbf{h}_\alpha - \left[ \mathbf{D}_\alpha \mathbf{N}^{(1)}(\mathbf{u}^0, \boldsymbol{\alpha}^0) \boldsymbol{\psi}^0 \right] \mathbf{h}_\alpha \end{pmatrix} \right\rangle_{u\psi} \quad (42)$$
$$- \left\{ \hat{P}_2(\mathbf{u}^0, \boldsymbol{\alpha}^0; \boldsymbol{\psi}_{i1}^{(2)}, \boldsymbol{\psi}_{i2}^{(2)}; \mathbf{h}_\alpha) \right\}_{\partial \Omega_x}.$$

In terms of the adjoint functions $\psi_{i1}^{(2)}$ and $\psi_{i2}^{(2)}$, the complete expression of the second-order mixed sensitivities $\delta S_i = (\delta S_i)_{direct} + (\delta S_i)_{indirect}$ is obtained by adding the above expression to the previously computed "direct effect term" from Eq. (26), to obtain

$$\delta S_i(\mathbf{u}^0, \boldsymbol{\alpha}^0, \boldsymbol{\psi}^0, \boldsymbol{\psi}_{i1}^{(2)}, \boldsymbol{\psi}_{i2}^{(2)}; \mathbf{h}_\alpha) = (\delta S_i)_{direct} + (\delta S_i)_{indirect}$$
$$= \frac{\partial}{\partial \alpha_i} \left\{ \left\langle \mathbf{D}_\alpha R(\mathbf{u}, \boldsymbol{\alpha}) - \mathbf{D}_\alpha \hat{P}(\mathbf{u}, \boldsymbol{\psi}, \boldsymbol{\alpha}) + \left\langle \boldsymbol{\psi}, \mathbf{D}_\alpha \mathbf{Q}(\boldsymbol{\alpha}) - \mathbf{D}_\alpha [\mathbf{L}(\boldsymbol{\alpha}) \mathbf{u}] \right\rangle_\psi, \mathbf{h}_\alpha \right\rangle_\alpha \right\}_{(\mathbf{u}^0, \boldsymbol{\psi}^0; \boldsymbol{\alpha}^0)} \quad (43)$$
$$+ \left\langle \left( \boldsymbol{\psi}_{i1}^{(2)\dagger}, \boldsymbol{\psi}_{i2}^{(2)\dagger} \right) \begin{pmatrix} \left\langle \mathbf{D}_\alpha \mathbf{Q}(\boldsymbol{\alpha}^0) - \mathbf{D}_\alpha [\mathbf{L}(\boldsymbol{\alpha}^0) \mathbf{u}^0], \mathbf{h}_\alpha \right\rangle_\alpha \\ \left\langle \mathbf{D}_{\alpha u}^2 R(e^0) - \mathbf{L}^*[(\boldsymbol{\alpha}^0) \boldsymbol{\psi}^0], \mathbf{h}_\alpha \right\rangle_\alpha \end{pmatrix} \right\rangle_{u\psi} - \hat{P}_2(\mathbf{u}^0, \boldsymbol{\alpha}^0; \boldsymbol{\psi}_{i1}^{(2)}, \boldsymbol{\psi}_{i2}^{(2)}; \mathbf{h}_\alpha).$$

It is convenient to call Eqs. (38) and (41) the *$2^{nd}$-level adjoint sensitivity system ($2^{nd}$-LASS) for the second-level adjoint function* $\boldsymbol{\psi}_i^{(2)} \equiv \left( \boldsymbol{\psi}_{i1}^{(2)}, \boldsymbol{\psi}_{i2}^{(2)} \right)$. *Note that the $2^{nd}$-LASS equations are independent of parameter variations* $\mathbf{h}_\alpha$. Evidently, *a single computation of the $2^{nd}$-LASS,* which yields the $2^{nd}$-level adjoint function $\boldsymbol{\psi}_i^{(2)} \equiv \left( \boldsymbol{\psi}_{i1}^{(2)}, \boldsymbol{\psi}_{i2}^{(2)} \right)$ for computing the corresponding the partial G-differential $\delta S_i(\mathbf{u}, \boldsymbol{\psi}, \boldsymbol{\alpha})$, *suffices to obtain subsequently, via* *inexpensive quadratures for computing the inner products shown in Eq. (43), the $i^{th}$ row/column* $\left\{\partial^2 R(\mathbf{u},\boldsymbol{\alpha})/\partial\alpha_j\partial\alpha_i\right\}_{(\mathbf{u}^0,\boldsymbol{\alpha}^0)}$, $j=1,...,N_\alpha$, *of the matrix of second-order response sensitivities.* For each index $i=1,...,N_\alpha$, *the second-level adjoint function* $\boldsymbol{\psi}_i^{(2)} \equiv \left(\boldsymbol{\psi}_{i1}^{(2)}, \boldsymbol{\psi}_{i2}^{(2)}\right)$ is the solution of the corresponding $2^{nd}$-LASS. Each $2^{nd}$-LASS comprises the same operators on its left-side; only the source terms on the rights sides of the $2^{nd}$-LASS differ from one another, since the right-sides correspond to the distinct having a distinct right-hand side (source) which stems from the partial G-differential $\delta S_i(\mathbf{u},\boldsymbol{\psi},\boldsymbol{\alpha})$, $i=1,...,N_\alpha$. Thus, *the exact computation of all second-order sensitivities,* $\left\{\partial^2 R(\mathbf{u},\boldsymbol{\alpha})/\partial\alpha_j\partial\alpha_i\right\}_{(\mathbf{u}^0,\boldsymbol{\alpha}^0)}$, *for* $i,j=1,...,N_\alpha$, *using the $2^{nd}$-ASAM requires at most* $N_\alpha$ *large-scale computations using the $2^{nd}$-LASS*, rather than $O(N_\alpha^2)$ large-scale computations as would be required by forward methods. It is also important to note that the construction and solution of the $2^{nd}$-LASS requires very little effort beyond that already invested in solving the original forward Eq.(1) for the state variable $\mathbf{u}(\mathbf{x})$ and the $1^{st}$-LASS, cf. Eq. (16) and (19), for computing the first-level adjoint function $\boldsymbol{\psi}$. This is because the left-sides of Eq. (41) comprise operators that are very similar to the left-side of the $1^{st}$-LASS. However, the right-sides of Eq. (41) comprise "source terms" that differ form the right-sides ("source terms") of the $1^{st}$-LASS. Notably, due to the symmetry properties of the second-order sensitivities $\left\{\partial^2 R(\mathbf{u},\boldsymbol{\alpha})/\partial\alpha_j\partial\alpha_i\right\}_{(\mathbf{u}^0,\boldsymbol{\alpha}^0)}$, $j=1,...,N_\alpha$, , the $2^{nd}$- ASAM provides an automatic, inherent, "solution verification" of the correctness and accuracy of the $2^{nd}$-level adjoint functions $\boldsymbol{\psi}_i^{(2)} \equiv \left(\boldsymbol{\psi}_{i1}^{(2)}, \boldsymbol{\psi}_{i2}^{(2)}\right)$ needed for the efficient and exact computation of the second-order sensitivities.

In the case of linear systems, i.e., *when Eqs. (1) and (2) are linear in* $\mathbf{u}(\mathbf{x})$, *so that* $\mathbf{N}[\mathbf{u}(\mathbf{x}),\boldsymbol{\alpha}(\mathbf{x})] \equiv \mathbf{L}[\boldsymbol{\alpha}(\mathbf{x})]\mathbf{u}(\mathbf{x})$ *and* $\mathbf{B}[\mathbf{u}(\mathbf{x}),\boldsymbol{\alpha}(\mathbf{x})] \equiv \mathbf{B}[\boldsymbol{\alpha}(\mathbf{x})]\mathbf{u}(\mathbf{x})$, *then the $2^{nd}$-LASS, i.e., Eqs. (41) and (38) will reduce to the corresponding equations for linear systems* [see Eqs. (38b) and (39), respectively, in Cacuci (2015)]. For linear systems, therefore, the solution of the $2^{nd}$-LASS simplifies considerably since, as discussed by Cacuci (2015):

(i) The left-side of the equation for determining the second-level adjoint function $\boldsymbol{\psi}_{i2}^{(2)}(\mathbf{x})$ will become the same as the left-side of the equation for determining $\mathbf{u}(\mathbf{x})$, except for a different source term; and

(ii) The left-sides of the equations to be solved for determining $\psi_{i1}^{(2)}(\mathbf{x})$ will become the same as the left-sides of the equations for determining the first-level adjoint function $\psi$, except (again) for a different source term.

## 4. CONCLUSIONS

This work has presented the "**S**econd-**O**rder **A**djoint **S**ensitivity **A**nalysis **M**ethodology ($2^{nd}$-*ASAM*) for computing exactly and efficiently the second-order functional derivatives of system responses (i.e., "system performance parameters") to the system's model parameters. The definition of "system parameters" used in this work include, in the most comprehensive sense, all computational input data, correlations, initial and/or boundary conditions, etc. The $2^{nd}$-*ASAM* builds on the "first-order adjoint sensitivity analysis methodology" ($1^{st}$-*ASAM*) for nonlinear systems originally introduced ([11, 12]) and developed ([13]) by Cacuci; see also Refs. [14-16]. It also extents the work in [1] to the consideration of not only nonlinear responses, but to fully nonlinear systems. This work has shown that, for *one functional-type response* of interest produced by a physical system comprising $N_\alpha$ parameters and $N_r$ responses, the $2^{nd}$-*ASAM requires one large-scale computation using the first-level adjoint sensitivity system ($1^{st}$-LASS) for obtaining all of the first-order sensitivities* $\{\partial R(\mathbf{u},\boldsymbol{\alpha})/\partial \alpha_i\}_{(\mathbf{u}^0,\boldsymbol{\alpha}^0)}$, $i=1,...,N_\alpha$, *followed by at most $N_\alpha$ large-scale computations using second-level adjoint sensitivity systems ($2^{nd}$-LASS) for obtaining exactly all of the second-order sensitivities* $\{\partial^2 R(\mathbf{u},\boldsymbol{\alpha})/\partial \alpha_j \partial \alpha_i\}_{(\mathbf{u}^0,\boldsymbol{\alpha}^0)}$, $i,j=1,...,N_\alpha$. In contradistinction, forward methods would require $(N_\alpha^2/2 + 3N_\alpha/2)$ large scale computations for obtaining all of the first- and second-order sensitivities, for all $N_r$ system responses, as follows: obtaining the first-order response sensitivities $\{\partial R(\mathbf{u},\boldsymbol{\alpha})/\partial \alpha_i\}_{(\mathbf{u}^0,\boldsymbol{\alpha}^0)}$, $i=1,...,N_\alpha$, requires $N_\alpha$ large-scale forward model computations, and obtaining the second-order sensitivities $\{\partial^2 R(\mathbf{u},\boldsymbol{\alpha})/\partial \alpha_j \partial \alpha_i\}_{(\mathbf{u}^0,\boldsymbol{\alpha}^0)}$, $i,j=1,...,N_\alpha$, requires an additional number of $N_\alpha(N_\alpha+1)/2$ large-scale computations.

Due to the symmetry properties of the second-order sensitivities $\left\{\partial^2 R(\mathbf{u},\boldsymbol{\alpha})/\partial\alpha_j\partial\alpha_i\right\}_{(\mathbf{u}^0,\boldsymbol{\alpha}^0)}$, $j=1,...,N_\alpha$, the $2^{nd}$-ASAM provides an automatic, inherent, "solution verification" of the correctness and accuracy of the $2^{nd}$-level adjoint functions $\boldsymbol{\psi}_i^{(2)} \equiv \left(\boldsymbol{\psi}_{i1}^{(2)},\boldsymbol{\psi}_{i2}^{(2)}\right)$ needed for the efficient and exact computation of the second-order sensitivities.

The *$2^{nd}$-ASAM for nonlinear systems* presented in this work enables the exact computation of all of the second-order response sensitivities (i.e., functional Gateaux-derivatives) to the large-number of parameters characteristic of the large-scale nonlinear systems of practical interest. The accompanying PART II [25] presents an illustrative application of the *$2^{nd}$-ASAM for nonlinear systems* to a paradigm nonlinear heat conduction problem that admits a unique analytical solution, thereby making transparent the mathematical derivations presented in this paper. Very importantly, this illustrative application will show that:

(i) The construction and solutions of the *second-level adjoint sensitivity systems ($2^{nd}$-LASS)* require very little additional effort beyond the construction of the first-level adjoint sensitivity system (*$1^{st}$-LASS*) needed for computing the first-order sensitivities; and

(ii) The symmetry properties of the second-order sensitivities $\left\{\partial^2 R(\mathbf{u},\boldsymbol{\alpha})/\partial\alpha_j\partial\alpha_i\right\}_{(\mathbf{u}^0,\boldsymbol{\alpha}^0)}$, $j=1,...,N_\alpha$, provide an automatic, inherent, "solution verification" of the correctness and accuracy of the $2^{nd}$-level adjoint functions $\boldsymbol{\psi}_i^{(2)} \equiv \left(\boldsymbol{\psi}_{i1}^{(2)},\boldsymbol{\psi}_{i2}^{(2)}\right)$ needed for the efficient and exact computation of the second-order sensitivities.

Although the *$2^{nd}$-ASAM for nonlinear systems* developed in this work was set in *real* (as opposed to complex) Hilbert spaces, since real Hilbert spaces provide the natural mathematical setting for computational purposes, This setting does *not* restrict, in any way, the generality of the concepts presented here. The *$2^{nd}$-ASAM for nonlinear systems* can be readily set in complex Hilbert spaces by simply changing some terminology, without affecting its substance. The use of the *$2^{nd}$-ASAM* to compute exactly all of the second-order response sensitivities to model input parameters is expected to enable significant advances in related scientific disciplines, particularly the areas of uncertainty quantification and predictive modeling, including model validation, reduced-order modeling, data assimilation, model calibration and extrapolation.

## ACKNOWLEDGMENTS

This work has been partially supported by Gen4Energy, Inc., under Contract 15540-FC51 with the University of South Carolina.